\documentclass[12pt]{article}
\usepackage{graphicx} 

\usepackage{amsmath}
\usepackage{amssymb}
\usepackage{amsthm}
\usepackage{latexsym}

\newcommand{\comment}[1]{}
\newtheorem{thm}{Theorem}        
\newtheorem{lemma}{Lemma}[section]
\newtheorem{remark}{Remark}[section]

\newtheorem{prop}{Proposition}[section]

\newtheorem{df}{Definition}[section]

\newtheorem{cor}{Corollary}[section]
\begin{document}

\date{}

\title{{\LARGE\sf Approach to Fixation for Zero-Temperature Stochastic Ising Models on the 
Hexagonal Lattice}}

\author{
{\bf Federico Camia}\\
{\small \tt federico.camia\,@\,physics.nyu.edu}\\
{\small \sl Department 
of Physics, New York University, New York, NY 10003, USA}\\
\and
{\bf Charles M.~Newman}\\
{\small \tt newman\,@\,courant.nyu.edu}\\
{\small \sl Courant Inst.~of Mathematical Sciences, 
New York University, New York, NY 10012, USA}\\
\and
{\bf Vladas Sidoravicius}\\
{\small \tt vladas\,@\,impa.br}\\
{\small \sl Instituto de Matematica Pura e Aplicada,
Rio de Janeiro, RJ, Brazil}\\
}

\maketitle

\begin{abstract} 
We investigate zero-temperature dynamics  on the hexagonal lattice 
${\Bbb H}$ for the homogeneous
ferromagnetic Ising model with zero external magnetic field and a disordered
ferromagnetic Ising model with a positive external magnetic field $h$.  
We consider both continuous time (asynchronous) processes and, in the
homogeneous case, also discrete time synchronous dynamics
(i.e., a deterministic cellular automaton), alternating
between two sublattices of ${\Bbb H}$.
The state space
consists of assignments of $-1$ or $+1$ to each site of ${\Bbb H}$,
and the processes
are zero-temperature limits of stochastic Ising ferromagnets with Glauber
dynamics and a random (i.i.d. Bernoulli) spin configuration at time $0$.  
We study the speed of convergence of the configuration $\sigma^t$ at time $t$
to its limit $\sigma^{\infty}$ and related issues.

\end{abstract}





\section{Introduction}

In this paper, we consider a number of continuous time Markov processes $\sigma^t, \; t\ge 0$,
with state space ${\cal S} = \{ -1, +1 \}^{\Bbb L}$ consisting of assignments of $-1$ or $+1$
to a regular lattice ${\Bbb L}$. In general, we take ${\Bbb L}$ to be 
$\Bbb H$, the hexagonal lattice in the 
plane, but occasionally we refer to results on $\Bbb Z^d$ or other lattices.
Later on, we will also consider a related discrete time process on $\Bbb H$.
A state $\sigma \in {\cal S}$ will also be called a spin configuration.
All the continuous time processes we consider are nearest-neighbor interacting particle systems
in which the dynamics may be constructed by means of independent, rate $1$, Poisson ``clock''
processes at  the sites $x$ of ${\Bbb L}$. If the Poisson clock at $x$ rings at time $t$, then a
\emph{spin flip} $(\sigma_x^{t+} = -\sigma_x^{t-})$ is considered; whether or not a flip actually
occurs depends on the values of $\sigma_x^{t-}$ and $\sigma_y^{t-}$
for each neighbor $y$ of $x$.

For the simplest model, the {\it homogeneous ferromagnet}, there is a flip with probability $1$ 
(or $1/2$ or $0$), if at $t-$, $\sigma_x^{\, \cdot}$ agrees with less than half (or exactly half or more
than half) of its neighbors. Note that the probability $1/2$ case does {\it not} occur in lattices like
$\Bbb H$ where every site has an {\it odd} number of neighbors so that there cannot be a tie among the
neighbors. We denote by $\omega$ a realization of the clocks rings (and tie-breaking coin tosses, 
if needed) and by $P_{\omega}$ the corresponding probability measure. In all our models we choose 
$\sigma^0$ according to a probability measure  $P_{\sigma^0}$ corresponding to i.i.d. ${\sigma^0_x}$'s
with $P_{\sigma^0}(\sigma_x^0=+1)=\lambda \in [0,1]$. In the Ising model context, $\lambda = 1/2$ is the
most important special case and corresponds to an initial ``quench from infinite temperature''.

In the physics literature on zero-temperature dynamics for the homogeneous ferromagnet with $\lambda = 1/2$,
an important quantity is the $P_{\sigma^0}\times P_{\omega}$ probability $p(t)$ that the origin has not
flipped at all by time $t$. In this situation, it has been known for a long time \cite{Arratia} that on 
$\Bbb Z$, almost surely (a.s.) every site flips infinitely often and hence $p(t) \to 0$ as $t \to + \infty$. 
More recently it was shown \cite{Derrida, DHP} that on $\Bbb Z$, $p(t)\sim t^{-3/8}$ and numerically 
\cite{DOS} that on $\Bbb Z^2$, $p(t) \sim t^{-\theta(2)}$ with $\theta (2) \approx 0.22$. 
This is consistent with the theorem \cite{NNS}
that also for $d=2$, a.s. every site flips infinitely often.
There is as yet no rigorous result for $\Bbb Z^d$ with $d>2$,
but numerical evidence \cite{Stauffer} suggests that for $d=3$, $\; p(t) \sim t^{-\theta (3)}$ with 
$0<\theta (3) < +\infty$ and raises at least the possibility that $p(t)\to p(\infty) >0$ for $d>4$. 
If the latter occurs, it may still be that $p(t) - p(\infty) \sim t^{-\theta (d)}$ with $0<\theta(d) < +\infty$. 
(We discuss the case $\lambda \neq 1/2$ below.)

The situation for lattices like $\Bbb H$ with no ties among neighbors is quite different. Here it is easy 
to see that $p(\infty)>0 \; $ (e.g., on $\Bbb H$, consider elementary hexagons with constant spin value at 
$t=0$ ) and a general theorem \cite{NNS} 
implies that a.s. every site flips only finitely many times so 
that $\sigma^{\infty} = \lim_{t\to \infty} \sigma^t$ exists. A natural question is then whether 
$p(t) - p(\infty)$ or its upper bound 

\begin{equation} \label{upperbound}
\tilde p (t) \equiv (P_{\sigma^0} \times P_{\omega}) (\text{origin flips after time $t$})
\end{equation}
tends to zero like $t^{-\phi}$ or exponentially fast or something
in between. 

Based on results for simpler quasi-one-dimensional 
lattices like $\Bbb Z \times \{0,1 \}$
\cite{NS2} and on results for the homogeneous tree of degree three \cite{Howard}, it was speculated in 
\cite{NS2} that on $\Bbb H$ convergence to the fixed, i.e., absorbing, state $\sigma^\infty$ may be
exponentially fast:
\begin{equation} \label{conjecture}
\tilde p(t) \le \exp \{- (ct)\}.
\end{equation}
This possibility is also supported by numerical evidence of Howard \cite{Howard,HN}.
Although this was originally considered
for $\lambda = 1/2$, it seems reasonable that such a conjecture should be valid 
for all $\lambda$. One of the main results of this 
paper is a proof of the following stretched exponential estimate with exponent $\gamma = 1/2$
for $\lambda$ not too close to 1/2 --- see Theorems \ref{exp_conv} and \ref{exp_conv1} in Section~2
below:
\begin{equation} \label{FSS}
\tilde p (t) \le \exp \{- (c't^\gamma) \}.
\end{equation} 
In this regard it is 
interesting to note that on $\Bbb Z^d$ with $d\ge 2$, and $\lambda$ 
sufficiently close to $1$ (respectively, $0$), 
it has been proved \cite{FSS} that $\sigma^t \to \sigma^\infty \equiv +1$
(resp., $\equiv -1$) with a similar stretched exponential upper bound, but
with $\gamma = \gamma (d) \in (0,1)$, and a similar lower bound for $d=2$.

Closely related to bounds like (\ref{conjecture}) or (\ref{FSS})
is the spatially localized mechanism for fixation, namely 
$\sigma_x^\cdot$ can no longer flip once it belongs to either a loop or ``barbell'' of constant
sign in $\Bbb H$. 
By loop we will always mean a simple loop (with no subloops).
A barbell consists of two disjoint loops connected by a path and we regard a loop as a degenerate
barbell. By studying the 
percolation properties of the final state $\sigma^\infty$ on the infinite lattice $\Bbb H$, it follows - 
see Prop. \ref{locally_fixated} - that for any $\lambda$, every site is in some
$\sigma^\infty$-barbell of constant sign.
Let $R$ denote the diameter of the smallest such barbell. We conjecture, that for any $\lambda$,

\begin{equation} \label{barbell}
(P_{\sigma_0} \times P_{\omega}) (R>r) \le \exp (-c'' r)
\end{equation}
for some $c''\in(0, \infty)$.
For $\lambda$ not close to 1/2, this
follows from the proofs of Theorems \ref{exp_conv} and \ref{exp_conv1} below. 
There would be an interesting application of the combination of (\ref{FSS}) and (\ref{barbell}) concerning
{\it overlap distributions} (cf. Sec. IX of \cite{NS3}), which we will briefly
describe at the end of this section after we introduce 
zero-temperature dynamics for disordered Ising models.

Our results for the continuous 
time dynamics for the homogeneous ferromagnet on ${\Bbb H}$ are
in Section 2. Then in Section 3, 
we analyze a discrete time dynamics on ${\Bbb H}$
(alternating between and synchronous within each of two sublattices)
that is the zero-temperature case of dynamics studied by Domany
\cite{domany} and that has also
been studied numerically by Nienhuis \cite{Nienhuis}.
This is really a deterministic cellular automaton and we obtain much stronger
results for it, such as exponential 
decay of $\tilde{p}$ for any $\lambda$, than we have
obtained for the usual continuous time dynamics.
We remark that some interesting results are obtained in
\cite{cns} on the continuum scaling limit of
the dependent percolation models generated by 
this cellular automaton when $\lambda = 1/2$.
Finally, in Section 4, we study 
continuous time dynamics for disordered Ising models on ${\Bbb H}$.

In disordered Ising models, the flip probability is determined by a (formal) Hamiltonian
\begin{equation} \label{Hamiltonian}
{\cal H} = - \sum_{ \langle x,y \rangle }
J_{x,y} \sigma_x \sigma_y - h \sum_x \sigma_x,
\end{equation}
where $\sum_{ \langle x,y \rangle }$ denotes the sum over all pairs of nearest
neighbor sites. 
We write ${\cal N}(x)$ for the set of nearest neighbors of $x$.
The homogeneous ferromagnet corresponds to the case where
\begin{itemize}
\item $J_{x,y} \equiv 1$ for all pairs of nearest neighbors $x,y$,
\end{itemize}
and we also take $h=0$ in that case, 
while disordered models correspond to the case where 
\begin{itemize}
\item the $J_{x,y}$'s are i.i.d. random variables.
\end{itemize}
We denote by $P_{\cal J}$ the distribution of
coupling realizations ${\cal J}$.
When the Poisson clock at $x$ rings, if
\begin{equation}
\Delta {\cal H}_x (\sigma) =
2 \sum_{y \in {\cal N}(x)} J_{x,y} \sigma_x \sigma_y + 2 h \sigma_x
\end{equation}
is negative,
then the flip is done with probability $1$;
if it is positive, then the flip is done with probability $0$.
If the energy change were zero, the flip would be done with 
probability $1/2$; however for the disordered models, we will
generally consider couplings whose common distribution is continuous,
so that there is zero probability for zero energy flips to occur. 
For such disordered models with continuous distributions of
couplings, it is also a consequence of a general result of \cite{NNS}
that $\sigma^{\infty} = \lim_{t \to \infty} \sigma^t $
exists. Again, we will be interested in the $t \to \infty$ behavior of 
$\tilde p (t)$, which here is the 
$P_{\sigma^0} \times P_{\omega} \times P_{\cal J}$
probability that the origin flips after time $t$.

Similar stochastic processes on different types of lattices have been studied
in various papers.
See, for example, \cite{CDN, FSS, GNS, NNS, NS2, NS1, NS} for models on ${\Bbb Z}^d$
and \cite{Howard} for a homogeneous ferromagnetic model on ${\Bbb T}_3$, the homogeneous
tree of degree three.
Such models are also discussed extensively in the physics literature, including
homogeneous and disordered ferromagnets as well as spin glasses, usually on the
${\Bbb Z}^d$ lattice.

Exponential decay of $\tilde{p} (t)$ has been obtained by Newman and Stein
\cite{NS2} for the homogeneous ferromagnet on the ``ladder'' lattice
(${\Bbb Z} \times \{ 0,1 \})$ and for continuously disordered models on ${\Bbb Z}$.
In both settings, the exponential decay is a result of a ``blocking'' condition present
at time zero and due to the initial random spin configuration in the case of the ladder
and to the structure of the random couplings in the case of the disordered ferromagnets
on ${\Bbb Z}$.
This blocking conditions break up the system into connected components that can have no
influence on each other, making the dynamics fundamentally local in nature.
The situation of the homogeneous ferromagnet on ${\Bbb T}_3$ is substantially different.
It is shown in \cite{Howard} that, if the density of $+1$ spins at time $0$ is large
enough, the system's agreement-inducing dynamics \emph{produces} enough ``fixated spins''
to break up the system and ensure exponential decay of $\tilde{p} (t)$.
The blocking mechanism in this case is dynamical.

We now turn to a discussion of overlap distributions for
disordered (and homogeneous) models. 
Let us denote by $\mu^\infty$ the probability distribution of the random final state $\sigma^\infty$ 
(induced by  $P_{\sigma_0} \times P_{\omega}$ or by
$P_{\sigma_0} \times P_{\omega} \times P_{\cal J}$ in disordered models) 
and let ${\sigma^{\infty}}' $ be a ``replica'' of $\sigma^\infty$
--- i.e., let the pair $(\sigma^\infty, {\sigma^{\infty}}' )$ 
 be distributed by  $ \mu^\infty\times \mu^\infty$. 
Then, by the spatial ergodic theorem, it follows \cite{NS3}
that the overlap random variable

\begin{equation} \label{overlap}
Q \equiv \lim_{L\to \infty} \frac{1}{|\Lambda_L \cap \Bbb H|} 
\sum_{x\in \Lambda_L \cap \Bbb H} \sigma_x^\infty {\sigma_x^{\infty}}',
\end{equation}
where $\Lambda_L$ denotes the $L\times L$ square centered at the origin
and $|\,\cdot\,|$ denotes cardinality, 
is a constant $g=g(\lambda)$ (with $g(1/2) =0$).
It is however a priori possible (see \cite{NS3} for a discussion) that for spin glass models,
where the common distribution of the couplings is symmetric about zero,
this triviality of the dynamical 
overlap distribution would not be so if one instead considered
the limit (in distribution) as $L\to \infty$ of the finite volume overlap

\begin{equation} 
Q_L \equiv  \frac{1}{|\Lambda_L \cap \Bbb H|} 
\sum_{x\in \Lambda_L \cap \Bbb H} \sigma_x^\infty [\Lambda_L] 
{\sigma_x^{\infty}}' [\Lambda_L] ,
\end{equation}
where  $\sigma^\infty [\Lambda_L]$ is the final state of the zero-temperature
Markov process restricted to the {\it finite volume} $\Lambda_L$ (with some boundary condition).
But it was conjectured in \cite{NS3} that $Q_L$ should rather
converge in probability to the same constant 
$g(\lambda)$. To prove this would require that $\sigma^\infty [\Lambda_L] \to
\sigma^\infty $ as $L\to \infty$ in such a way 
that $\sigma^\infty [\Lambda_L]$ and $\sigma^\infty $ only disagree 
in a small ``boundary layer'' (of area $o (L^2)$). 
In Corollary \ref{interchange} we do show that $\sigma^\infty [\Lambda_L] \to
\sigma^\infty $, and the extra uniformity needed would 
follow easily by combining conjectures (\ref{FSS}) and (\ref{barbell}).
We have in fact this result for $\lambda$ not  too close to 1/2, but 
the physically more interesting case is $\lambda = 1/2$.  
Of course, the most interesting result would be for a spin glass
rather than for a ferromagnet, homogeneous or otherwise.

\bigskip
\bigskip

\section{Homogeneous ferromagnet}

In this section, we consider the homogeneous ferromagnet on $\Bbb H$ with zero external magnetic
field, that is, $J_{x,y} \equiv 1$ for all pairs of nearest neighbors
$x,y$ and $h=0$.
Let $\lambda$ be the density of $+1$ spins in $\sigma^0$.

Our first result concerns the almost sure absence of percolation of both
$+1$ and $-1$ spins in the final configuration $\sigma^{\infty}$ (or in
$\sigma^t$) when $\lambda = 1/2$.
We note that it has been proved \cite{HN} that in $\sigma^{\infty}$, the mean cluster
size is infinite (see also Theorem \ref{mean_c_s} below).

\begin{prop} \label{percolation}
In the homogeneous ferromagnet,
if $\lambda = 1/2$ and $t \in [0, \infty]$, then for almost every $\sigma^0$ and $\omega$
there is no percolation in $\sigma^t$ of either $+1$ or $-1$ spins.
\end{prop}

{\bf Proof.} First note that the measure $\mu^t$ describing the state
$\sigma^t$ of the system at time $t \in [0, \infty]$ is invariant and
ergodic under any ${\Bbb H}$-automorphism. 
This is so  because the same is true for both
$P_{\sigma^0}$ and $P_{\omega}$ and hence also for $P_{\sigma^0} \times P_{\omega}$.
Applying a result of Harris \cite{Harris, Liggett}, we also have that
$\mu^t$ satisfies the FKG property, i.e., 
increasing functions of the spin
variables are positively correlated (this follows from the FKG property of
$P_{\sigma^0}$ and the attractiveness of the Markov process).
Then it follows from a result of Gandolfi, Keane and Russo (stated in \cite{GKR}
for ${\Bbb Z}^2$, but also valid for ${\Bbb H}$) that if 
percolation of, say, $+1$ sites were to occur, all the $-1$ clusters 
would have to be finite. 
Because of the symmetry of the model under a global spin flip, however,
percolation of $+1$ sites with positive probability implies the same for 
$-1$ sites. 
Then, using the ergodicity of the measure, we would see simultaneous
percolation of both signs, thus obtaining a contradiction. \fbox{} \\

\begin{remark}
Numerical evidence for the homogeneous ferromagnet \cite{HN} suggests that
there is plus (resp., minus) percolation in $\sigma^{\infty}$ for all
$\lambda>1/2$ (resp., $\lambda<1/2$).
For such a result in the case of \emph{synchronous} dynamics, see Prop. \ref{perc} below.
\end{remark}

\begin{remark} 
The proof of Proposition \ref{percolation} works in a more general context
and was used in \cite{CDN} to get a similar result for the homogeneous ferromagnet
on ${\Bbb Z}^2$ (see Proposition 3.2 of \cite{CDN}).
The proof also shows that the symmetric Bernoulli product measure may be
replaced by any distribution for $\sigma^0$ which (a) is symmetric under
$\sigma^0 \rightarrow - \sigma^0$, (b) is invariant and ergodic under any (nontrivial)
${\Bbb H}$-automorphism, and (c) satisfies the FKG property.
\end{remark}

Let us call a configuration $\sigma \in {\cal S}$ \emph{locally fixated} if
for each $x \in {\Bbb H}$, there exists a (finite) lattice animal $A_x$ containing $x$
such that for any $\tilde{\sigma} \in {\cal S}$ that coincides with $\sigma$ on
$A_x$, $\Delta {\cal H}_y (\tilde{\sigma}) > 0$ for every $y \in A_x$, i.e., $y$
agrees with a strict majority of its neighbors.

\begin{prop} \label{locally_fixated} In the homogeneous ferromagnet,
for any $\lambda$, for almost every $\sigma^0$ and $\omega$, the final configuration 
$\sigma^{\infty}$ is locally fixated.
\end{prop}

{\bf Proof.} Let's first consider the case $\lambda = 1/2$.
In this case the claim follows immediately from Proposition \ref{percolation} for
$t = \infty$ and the fact that $\sigma^{\infty}$ is a.s. fixated.
In fact, for each $x \in {\Bbb H}$, it is enough to take $A_x$ to be the almost 
surely finite (e.g. plus) cluster at $x$.
We now write $\sigma^{\infty} (\lambda)$ to indicate dependence on the parameter $\lambda$
and consider $\lambda > 1/2$. 
Let's couple $\sigma^t (\lambda)$ with $\sigma^t (1/2)$, where $\sigma^t (\lambda)$
and $\sigma^t (1/2)$ have the same dynamics realization $\omega$ and 
$\sigma^0 (1/2) \leq \sigma^0 (\lambda)$ where $\leq$ indicates the natural partial
order.
From the attractiveness of the dynamics it follows that $\sigma^t (1/2) \leq \sigma^t (\lambda)$
for all $t \in [0, \infty]$.
Hence, each site $x$ such that (both $\sigma_x^{\infty} (\lambda) = +1$ and)
$\sigma_x^{\infty} (1/2) = +1$ is contained in an almost surely finite cluster
$C'_x$ with $C'_x \subset C_x$, where $C_x$ is the cluster at $x$ in
$\sigma^{\infty} (\lambda)$ and $C'_x$ is the cluster at $x$ in $\sigma^{\infty} (1/2)$.
Then, we let $A_x = C'_x$.
For a site $x$ such that $\sigma_x^{\infty} (1/2) = -1$, the cluster $C'_x$ that
contains $x$ in $\sigma_x^{\infty} (1/2)$ is almost surely finite and therefore
surrounded by $C^+_x$, the finite union of all neighboring plus clusters.
(This union is actually a single (connected) cluster, as can be seen using arguments
connected to those used in the proof of Theorem \ref{mean_c_s}, but we will not use that fact.)
Since $\sigma^{\infty} (1/2) \leq \sigma^{\infty} (\lambda)$, $C^+_x$ is still plus in 
$\sigma^{\infty} (\lambda)$. 
We let $A_x = C'_x \cup C^+_x$.
This completes the proof for $\lambda > 1/2$.
The case $\lambda < 1/2$ follows by symmetry. \fbox{} \\

The property of $\sigma^{\infty}$ of being locally fixated means that the dynamics 
produces finite clusters that are stable for local reasons. 
This, in turn, has the following interesting consequence (see, for example,
\cite{NS3} for some physics motivation).
Consider the sequence $\Lambda_L$ of squares and denote by $\sigma^t [\Lambda_L]$
the Markov process defined on the sublattice $\Lambda_L \cap {\Bbb H}$ (with some
boundary condition). 
Then it is easy to see that 
for every $x \in {\Bbb H}$,
and using a natural coupling between the processes for different $L$'s (see the proof
of Corollary \ref{interchange} below), the following is true almost surely:
\begin{equation} \label{limits}
\lim_{t \to \infty} \left( \lim_{L \to \infty} \sigma^t_x [\Lambda_L] \right)
= \lim_{t \to \infty} \sigma^t_x = \sigma^{\infty}_x.
\end{equation}
This follows from a finite speed of propagation of information argument
(an argument of this kind goes back to Harris \cite{Harris1} and is necessary even
to prove that the stochastic process itself is well defined; for a reference, see
\cite{Harris1} or p. 119 of \cite{Durrett}).
It is not clear though whether the two limits ($t \to \infty$ and $L \to \infty$)
commute.
In fact, on a homogeneous tree, interchanging the space and time limits can produce
a different result (for example, choosing the plus boundary condition).
But in our case, Proposition \ref{locally_fixated} implies the following corollary:

\begin{cor} \label{interchange} Almost surely,
\begin{equation}
\lim_{L \to \infty} \left( \lim_{t \to \infty} \sigma^t_x [\Lambda_L] \right)
= \lim_{t \to \infty} \left( \lim_{L \to \infty} \sigma^t_x [\Lambda_L] \right)
= \sigma^{\infty}_x.
\end{equation}
\end{cor}

{\bf Proof}. For almost every initial configuration $\sigma^0$ and realization of the
dynamics $\omega$, $\sigma^{\infty} (\sigma^0, \omega)$ is locally fixated and
therefore, at time $t = \infty$, $x$ belongs to a lattice animal $A_x$ which is stable
for local reasons.
Call $T_{A_x}$ the time it takes for all the sites $y \in A_x$ to fixate.
$T_{A_x} < \infty$ almost surely.
Let us now introduce a coupling between $\sigma^t [\Lambda_L]$ for every $L$ and
$\sigma^t$ constructed as follows:
1) the dynamics realizations for $\sigma^t [\Lambda_L]$ and $\sigma^t$ coincide in
$\Lambda_L \cap {\Bbb H}$, and 2) $\sigma^0_x = \sigma^0_x [\Lambda_L] ~ \forall x \in
\Lambda_L \cap {\Bbb H}$. 
If, for any given $\varepsilon > 0$, one can take $L$ large enough so that
$A_x \subset \Lambda_L \cap {\Bbb H}$ with probability at least $1 - \varepsilon$ and so that moreover
\begin{equation}
P( \sigma^{T_{A_x}}_{A_x} [\Lambda_L] \neq \sigma^{T_{A_x}}_{A_x} ~ | ~
\sigma^0 [\Lambda_L] = \sigma^0_{\Lambda_L} ) < \varepsilon,
\end{equation}
where $\sigma_B$ is the configuration $\sigma \in {\cal S}$ restricted to the set $B$,
the corollary follows.
But for any $\varepsilon > 0$, such an $L = L(\varepsilon)$ exists because of the finiteness of $T_{A_x}$
and the finite speed of propagation of information (once again, see \cite{Harris1} or p. 119 of \cite{Durrett}),
and this concludes the proof.  \fbox{} \\

Before we state our next result, we need some notation.
Call any minimal hexagon connecting six sites of ${\Bbb H}$ a \emph{cell}.
Each site is contained in three cells.
Given any connected subset $A$ of ${\Bbb H}$, the \emph{internal energy} of $A$
(denote by $e_A$) is defined as
\begin{equation}
e_A = - \sum_{ \stackrel{\langle x,y \rangle}{x,y \in A} } J_{x,y} \sigma_x \sigma_y,
\end{equation}
the sum being performed over pairs of neighboring sites, each pair counted once.

\begin{prop} \label{max}
In the homogeneous ferromagnet,
starting from any initial state (i.e., spin configuration), each site flips at most
$8$ times.
\end{prop}

{\bf Proof.} By the translation invariance of the model, it is enough to prove
the claim for the origin $0$.
Let $A$ be the union of all sites belonging to the three cells that contain the origin. Then $-15 \leq e_A \leq 15$.
Every spin flip of the origin or of one of its three neighbors lowers by at least two
units the internal energy of $A$, while a spin flip of any one of the other sites in $A$
can lower the internal energy $e_A$ or leave it unchanged, but can never raise it.
Therefore $e_A$ can never increase and moreover it strictly decreases by at least two
units whenever the origin or one of its neighbors flips.
Notice that after the origin has flipped the first time, every other spin flip of the origin has 
to be preceded by the spin flip of at least one of its neighbors. 
Thus every spin flip of the origin after the first one corresponds to a decrease in
$e_A$ of at least four units.
Given that the maximum change in $e_A$ (the energy when all edges are
unsatisfied minus the energy when all edges are satisfied) is $30$, the origin can
flip at most $8$ times.
%

\begin{remark}
An analogous result holds true for the homogeneous ferromagnet on the ladder 
lattice $({\Bbb Z} \times \{ 0,1 \})$. The proof is the same once the
corresponding cell has been properly defined.
\end{remark}

Recall that $\tilde{p}(t)$ denotes the probability that the origin flips after time $t$.
The following result says that, starting from a large enough density $\lambda$ of
$+1$ sites, the time it takes for $\sigma^t$ to converge to $\sigma^{\infty}$
has a tail decreasing at least as fast as a stretched-exponential.
We denote by $p_c^{site}$ the critical value for independent site percolation
on ${\Bbb H}$.
The symbols $c$, $c'$ and $c''$ will denote generic positive constants, whose values
may be different in different parts of the paper.

\begin{thm} \label{exp_conv}
In the homogeneous ferromagnet, if $\lambda > p_c^{site}$, then 
\begin{equation}
\tilde{p}(t) \leq e^{- c \sqrt{t}}
\end{equation}
for some $c \in (0, \infty)$.
\end{thm}

{\bf Proof.} If the density $\lambda$ of $+1$ sites is larger than $p_c^{site}$, there
exists almost surely an infinite cluster of $+1$ sites at time zero.
Every doubly-infinite path or closed loop of sites of the same sign is stable
for the dynamics. 
By standard percolation arguments, there are loops of stable $+1$ sites that
break up the lattice into finite subsets, and the probability that such subsets
are large is stretched-exponentially small.
To be more precise, if $D_0$ is the subset that contains the origin and is
surrounded by the smallest such loop, then
\begin{equation} \label{exp}
P_{\sigma^0} (|D_0| > N) \leq e^{- c' \sqrt{N}}
\end{equation}
for some $c' \in (0, \infty)$.
To sketch the arguments that lead to this, let $R_L$ be any $L \times 3L$ rectangle
in ${\Bbb R}^2$ and consider the event $A(R_L)$ that the set
$R_L \cap {\Bbb H}$ contains a plus-crossing joining its two sides of length $L$.
By standard arguments (\cite{Grimmett}; see also \cite{Russo} and \cite{SW}),
\begin{equation} \label{crossing}
P_{\sigma^0} (A(R_L)) \geq 1 - e^{- \beta L}
\end{equation}
for $L$ large enough and some $\beta \in (0, \infty)$.
Let $B(L)$ be the set $([-L/2, L/2] \times [-L/2, L/2]) \cap {\Bbb H}$ and
$S(L)$ be the event that there is a plus loop in $B(3L) \setminus B(L)$ surrounding
the origin.
Then, writing $B(3L) \setminus B(L)$ as a union of four $L \times 3L$ rectangles
and using (\ref{crossing}), we have
\begin{eqnarray}
P_{\sigma^0}(|D_0| \geq (3L)^2) & \leq & 1 - P_{\sigma^0}(S(L)) \nonumber \\
 & \leq & 4 \, e^{- \beta L} \nonumber \\
 & \leq & e^{- c'' L}
\end{eqnarray}
for $L$ large enough and some $c'' \in (0, \infty)$, which yields (\ref{exp}).

By Proposition \ref{max}, the maximum number of spin flips allowed inside $D_0$ is
$8 \, |D_0|$, where $|D_0|$ is the cardinality of $D_0$.
Moreover, as long as a single unstable site is present, it will flip at its next
clock ring.
It follows that an upper bound for the time for the sites in $D_0$ to fixate is
given by a sum of $8 \, |D_0|$ independent exponential (mean one) random variables
$T_1, T_2, \dots$ . 
More precisely, for any $\alpha >0$,
\begin{eqnarray}
\tilde{p}(t) & \leq & P_{\sigma^0} (|D_0| > [\alpha \, t]/8) + 
\max_{1 \leq n \leq [\alpha \, t]} {\rm Prob}(T_1 + \ldots + T_n > t) \nonumber \\
& = & P_{\sigma^0}(|D_0| > [\alpha \, t]/8) + {\rm Prob}(T_1 + \ldots + T_{[\alpha \, t]} 
> t). \label{p_t} 
\end{eqnarray}
By choosing $0 < \alpha < 1$, the first term at the end of (\ref{p_t}) can be seen to
be exponentially small in $\sqrt{t}$, using (\ref{exp}), and the second term is even
smaller (exponentially small in $t$) by standard large deviation arguments. \fbox{} \\

The next result is a slight improvement of the previous theorem.

\begin{thm} \label{exp_conv1}
There exists $\bar{\lambda} < p_c^{site}$ such that
if $\lambda \in (\bar{\lambda}, 1]$ in the homogeneous ferromagnet, then
\begin{equation}
\tilde{p}(t) \leq e^{- c \sqrt{t}}
\end{equation}
for some $c \in (0, \infty)$.
\end{thm}

{\bf Proof.} To prove the theorem, we will use a general result of 
Aizenman and Grimmett \cite{AG, Grimmett} on enhanced percolation. 
To do this, first of all
let us partition the hexagonal lattice into two sublattices ${\cal A}$
and ${\cal B}$ (with the origin of ${\Bbb H}$ in ${\cal B}$) in such a way that
all three neighbors of any site in
${\cal A}$ (resp., ${\cal B}$) are in ${\cal B}$ (resp., ${\cal A}$).
Any two sites of either ${\cal A}$ or of ${\cal B}$ have no edge of the original
hexagonal lattice in common.
By joining two sites of ${\cal A}$ whenever they are next-nearest neighbors
in the hexagonal lattice (two steps away from each other), we get a triangular
lattice (the same with ${\cal B}$).
For each site $x \in {\cal A}$, write the exponential random variable 
representing the time of the first ring of the Poisson clock at $x$ as a sum of three
independent identically distributed random variables (this is possible
because the exponential distribution is infinitely divisible).
Assign each one of these three random variables to one of the three neighbors of site
$x$ (such neighbors belong to the sublattice ${\cal B}$).
We now construct our enhanced percolation process in the following way:
\begin{itemize}
\item A site in the sublattice ${\cal A}$ is \emph{open} if its spin is plus.
\item A site in the sublattice ${\cal B}$ is \emph{open} if 
\begin{enumerate}
\item its spin is plus, or 
\item the random variable representing the first clock ring of that site is smaller
than each of the three random variables it was assigned from its three neighbors
in ${\cal A}$ and smaller than one.
\end{enumerate}
\end{itemize}
This last condition implies that that site in ${\cal B}$ will make its first attempt at
flipping before any of its three neighbors in ${\cal A}$.
The enhancement only takes place at sites of the sublattice ${\cal B}$ and is activated
at site $x$ with strictly positive probability. 

The enhancement is easily seen to be \emph{essential}, as defined in \cite{AG},
and therefore the result of Aizenman and Grimmett \cite{AG} (see also Section 3.3 of
\cite{Grimmett}) can be applied.
Thus, the critical value $p_c^{enh}$ for the density of $+1$ sites in $\sigma^0$ to yield
percolation of open sites is strictly lower
than that for independent site percolation: $p_c^{enh} < p_c^{site}$.
To conclude the proof we need to show that $\lambda > p_c^{enh}$ implies
stretched exponential decay of $\tilde{p} (t)$.
This is so because we claim that any loop of open sites in ${\Bbb H}$
is a loop of plus sites at time one, from which stretched-exponential decay
with exponent $1/2$ follows by using essentially the same arguments already used in the
proof of Theorem \ref{exp_conv}.

To see why any loop of open sites at time zero becomes a loop of plus
sites at time one, notice that by time one each open site of the lattice ${\cal B}$
which was not plus at time zero will have attempted a spin flip before any of its
neighbors and therefore will have flipped to plus. \fbox{} \\

\begin{remark} We conjecture that $\tilde{p} (t)$ decays exponentially for all
values of the initial density of $+1$ spins $\lambda$, including $1/2$
(see also \cite{Howard,NS2}).
This is also supported by simulation results of Howard \cite{Howard,HN}.
\end{remark}

\bigskip
\bigskip

\section{Synchronous dynamics}

In this section we study a different kind of dynamics, in discrete time, for the homogeneous
ferromagnet on $\Bbb H$.
This is the zero-temperature case of a dynamics studied by Domany \cite{domany}.
The rules (but not the timing) for updating the spins are the same as in the previous
section and we will use the same notation. The hexagonal lattice is partitioned
into two sublattices ${\cal A}$ and ${\cal B}$ (as in the proof of Theorem \ref{exp_conv1})
in such a way that the set ${\cal N} (x)$ of all three neighbors of a site $x$ in
${\cal A}$ (resp., ${\cal B}$) is in ${\cal B}$ (resp., ${\cal A}$).
By joining two sites of ${\cal A}$ whenever they are next-nearest neighbors
in the hexagonal lattice (two steps away from each other), we get a triangular
lattice (the same with ${\cal B}$).
The synchronous dynamics is such that all the sites in the sublattice ${\cal A}$
(resp., ${\cal B}$) are updated simultaneously.
Since this is a discrete time dynamics, our stochastic process will be denoted by
$\sigma^n$, with $n \in \{ 0,1,2,\dots \}$; $\sigma^0$ is still chosen from the
Bernoulli distribution $P_{\sigma^0}$ with density $\lambda$ of $+1$ spins.
The first update $\sigma^1$ will be for the sublattice ${\cal A}$.

\begin{prop} \label{perc}
If $\lambda > 1/2$ (resp., $< 1/2$), there is percolation of $+1$ (resp., $-1$)
spins in $\sigma^n$ for any $n \in [1, \infty]$.
\end{prop}

{\bf Proof.} We will only give the proof for $\lambda > 1/2$, since the case
$\lambda < 1/2$ is the same by symmetry.
If $\lambda > 1/2$, since the critical value for independent (Bernoulli) percolation
on the triangular lattice is exactly $1/2$, there is at time zero percolation of
$+1$ spins in the two triangular sublattices ${\cal A}$ and ${\cal B}$.
Any site $x \in {\cal A}$ with two $+1$ neighbors in ${\cal B}$ will become $+1$
when it updates, but since
$+1$ spins percolate (and form doubly-infinite paths) in ${\cal B}$,
at time $0$, when the sites of ${\cal A}$ are updated at time $1$, doubly-infinite paths of $+1$ spins
will be created in the hexagonal lattice and these are then stable.
Therefore, for $n \geq 1$, there is percolation of $+1$ spins in the hexagonal lattice. \fbox{} \\
%

\begin{thm} \label{syn_exp}
For any $\lambda$, $\sigma^n$ converges to $\sigma^{\infty}$ exponentially fast
for the synchronous dynamics in the sense that
\begin{equation} \label{syn_expon}
\tilde{p}^{\cal A} (n), \tilde{p}^{\cal B} (n) \leq e^{-cn}
\end{equation}
for some $c \in (0, \infty)$, where $\tilde{p}^{\cal A} (n)$ is the probability that
a deterministic site in ${\cal A}$ flips after time $n$ and similarly for $\tilde{p}^{\cal B} (n)$.
\end{thm}

{\bf Proof.} Without loss of generality, we assume that $0 < \lambda < 1$.
First of all we make the following observations: 
\begin{itemize}
\item the values of the spins in the sublattice ${\cal A}$ at time zero are irrelevant,
since at time $1$, after the first update, those values are uniquely determined by the
values of the spins in the sublattice ${\cal B}$,
\item once the initial spin configuration in the sublattice ${\cal B}$ is chosen, the
dynamics is completely deterministic.
\end{itemize}

We now concentrate on the sublattice ${\cal B}$. 
It is easy to see that the deterministic dynamics in ${\cal B}$
(observed at even times $2n$) is of nearest neighbor type.
We also claim that each site $x \in {\cal B}$ can flip at most \emph{one} time.
To see this, assume without loss of generality that a $-1$ spin at $x \in {\cal B}$
flips to $+1$ at time $2n$ with $n \geq 1$.
For this to happen, two of the three ${\Bbb H}$-neighbors of $x$ ($y_1, y_2 \in {\cal A}$)
must be $+1$ at time $2n$.
Therefore, at time $2n-1$, $y_1$ and $y_2$ each need to have two $+1$ ${\Bbb H}$-neighbors
in ${\cal B}$.
This implies that at time $2n-1$ all the sites in the hexagon containing $x, y_1, y_2$,
except $x$, are $+1$.
When, at time $2n$, the spin at $x$ flips to $+1$, a stable loop is formed, and site $x$ fixates.

Let us consider a loop $\gamma$ in the triangular sublattice ${\cal B}$, written as an ordered
sequence of sites $(y_0, y_1, \dots, y_n)$ with $n \geq 3$, which are distinct except that $y_n=y_0$.
For $i=1, \dots, n$, let $\zeta_i$ be the unique site in ${\cal A}$ that is an ${\Bbb H}$-neighbor
of both $y_{i-1}$ and $y_i$.
We call $\gamma$ an \emph{s-loop} if $\zeta_1, \dots, \zeta_n$ are all distinct.
Similarly, a (site-self avoiding) path $(y_0, y_1, \dots, y_n)$ in ${\cal B}$ is an \emph{s-path} if
$\zeta_1, \dots, \zeta_n$ are all distinct. 
Notice that any path in ${\cal B}$ (seen as a collection of sites) contains an s-path.
We will abuse our terminology slightly and also call a doubly-infinite s-path an s-loop.
An s-loop of constant sign is stable for the dynamics since at the next update of ${\cal A}$
the presence of the constant sign s-loop in ${\cal B}$ will produce a stable loop of that sign
in the hexagonal lattice.
A triangular loop $x_1, x_2, x_3 \in {\cal B}$ with a common
${\Bbb H}$-neighbor $\zeta \in {\cal A}$ is called a \emph{star}; it is not an s-loop.
A triangular loop in ${\cal B}$ that is not a star is an s-loop and will be called an \emph{antistar},
while any loop in ${\cal B}$ that contains more than three sites contains an s-loop.
A connected cluster in ${\cal B}$
that does not contain any loops other than stars (i.e., it doesn't contain s-loops)
will be called an \emph{s-tree} (notice that this does \emph{not} correspond to the usual definition
of a tree).
An s-tree can contain more than one star, but any two distinct stars cannot have sites in common or
they would together contain an s-loop.

Consider a specified site $x$ (e.g., the origin) in the triangular sublattice ${\cal B}$ and denote
by $C^{\cal B}_x$ its ${\cal B}$-cluster of constant sign (at some time).
The site $x$ can be of three types:
\begin{itemize}
\item $x$ belongs to an s-loop or to an \emph{s-barbell} which consists of two s-loops of
constant sign connected by a path (and therefore by an s-path) of the same sign,
in which case it is called \emph{frozen},
\item $x$ has only one ${\cal B}$-neighbor of its sign or else exactly two ${\cal B}$-neighbors
$y_1, y_2$ of its sign and $x, y_1, y_2$ share a common ${\Bbb H}$-neighbor (so that $(x, y_1, y_2)$
is a star), in which case it is called \emph{hot},
\item $x$ is neither frozen nor hot, in which case it is called \emph{warm}.
\end{itemize}
Notice that:
\begin{itemize}
\item a frozen site will never flip at any later time; a hot site will flip at the next time step;
a warm site will not flip at the next attempt, but might flip at a subsequent time,
\item any cluster that is not completely frozen must contain at least one hot site
(because the removal of all frozen sites results in one or more s-trees),
\item a cluster can grow but every newly
added site will be frozen (as shown in the argument above
that each site can only flip once), each newly added site
being previously a hot site of a cluster of opposite sign.
\end{itemize}

To analyze when a site $x \in {\cal B}$ will flip for the last time, we consider its ${\cal B}$-cluster
at time zero.
According to the above considerations, if $x$ is hot at time zero, then it will fixate at the first
update of ${\cal B}$; if it has four or more ${\cal B}$-neighbors of its sign at time zero,
then it belongs to an s-loop and is frozen.
There are only three possibilities left to be considered:
\begin{enumerate}
\item $x$ has exactly two ${\cal B}$-neighbors, $y_1$ and $y_2$, of its own sign and such that $x,y_1,y_2$
don't form a star.
\item $x$ has exactly three ${\cal B}$-neighbors, $y_1$, $y_2$ and $y_3$, 
of its own sign and such that none
of the combinations $x,y_i,y_j$, for $i,j=1,2,3$, forms a star.
\item $x$ has exactly three ${\cal B}$-neighbors, $y_1$, $y_2$ and $y_3$,
of its own sign and such that
$x,y_1,y_2$ form a star and $y_3$ is not a ${\cal B}$-neighbor of $y_1$ or $y_2$.
\end{enumerate}

Now, let $U$ be a specific nonempty subset of the set ${\cal N}^{\cal B}(x)$ of all six
${\cal B}$-neighbors of $x$.
We define the ``partial cluster'' $C^{\cal B}_{x,U}$ to be the set of sites $y \in {\cal B}$ such that
there is a (site-self avoiding) path $x_0=x, x_1, \dots, x_n=y$ of constant spin value (at time zero)
in ${\cal B}$ with $x_1 \in U$; i.e. $C^{\cal B}_{x,U}$ is the set of sites in $C^{\cal B}_x$ that can
be reached from $x$ by paths that start off by going from $x$ to a neighbor in $U$.
Note that if we are considering one of the cases 1, 2 or 3 described above, we can define a
\emph{branch} of $C^{\cal B}_x$ as $C^{\cal B}_{x,U}$ with $U$ either a singleton or, in case 3, 
also the doubleton $\{ y_1,y_2 \}$. In case 1, $C^{\cal B}_x$ has two branches; in case 2, it has three branches;
in case 3, it has two branches, one of which has a doubleton $U$.
For $x$ not to be frozen, the branches of $C^{\cal B}_x$ must be distinct (otherwise $x$ would be in
an s-loop), and all but one of them must be s-trees, i.e., must contain no loops other than stars
(otherwise $x$ would be in an s-barbell).
In this case, $x$ may eventually flip, and the time at which it will do so is bounded by the length of
the longest s-path contained in one of the branches of $C^{\cal B}_x$ that is an s-tree at time zero.

To complete the proof, it suffices to show that, for any fixed $U$, there is some $\beta>0$ and 
$K < \infty$ such that
\begin{equation} \label{antistar}
P(|C^{\cal B}_{x,U}| \geq n \text{ and $C^{\cal B}_{x,U}$ contains no antistar}) \leq K \, e^{- \beta n}.
\end{equation}
To prove (\ref{antistar}), we partition ${\cal B}$ into disjoint antistars and denote by $\tau$
the collection of these antistars.
We do an algorithmic construction of $C^{\cal B}_{x,U}$ (as in, e.g., \cite{FN}), where the order of
checking the sign of sites is such that when the first site in an antistar from $\tau$ is checked
(and found to have the same sign as $x$), then the other two sites of that antistar are checked next.
Without loss of generality, we assume that $\sigma^0_x = +1$.
Then standard arguments show that the probability in (\ref{antistar}) is bounded by
$K \, (1 - \lambda^3)^{(n/3)}$, which yields (\ref{syn_expon}) for $\tilde{p}^{\cal B} (n)$.

If a site $z$ belongs to the triangular sublattice ${\cal A}$, it is clear that once its three
${\Bbb H}$-neighbors have fixated, $z$ will certainly be fixated at the next update of ${\cal A}$,
and therefore the exponential bound (\ref{syn_expon}) holds also for $\tilde{p}^{\cal A} (n)$. \fbox \\

\begin{thm} \label{mean_c_s}
If $\lambda = 1/2$ in the synchronous dynamics \emph{:} 
\begin{enumerate}
\item For $n \in [0, \infty]$, there is no percolation in $\sigma^n$ of either $+1$ or
$-1$ spins, for almost every $\sigma^0$.
\item There is infinite mean cluster size in $\sigma^n$ for any $n \in [1, \infty]$ \emph{:}
for any $x \in {\Bbb H}$,
\begin{equation} \label{inf_mean}
E(|C_x(n)|) = \infty.
\end{equation}
\end{enumerate}
\end{thm}

{\bf Proof.} The proof of the first claim is the same as the proof of Proposition
\ref{percolation}.
%
To prove the second claim, let us first notice that $C_x(n) \cap {\cal B} = C_x^{\cal B}(n)$
for any $x \in {\cal B}$ and odd $n \geq 1$ and
$C_x(n) \cap {\cal A} = C_x^{\cal A}(n)$
for any $x \in {\cal A}$ and even $n \geq 2$.
Thus percolation in ${\cal B}$ for odd $n \geq 1$ (resp., in ${\cal A}$ for even $n \geq 2$)
would imply percolation
in the hexagonal lattice at the same $n$.
Therefore, by the first part of the theorem, at any odd $n \geq 1$
(resp., even $n \geq 2$) there is no
percolation in ${\cal B}$ (resp., ${\cal A}$).
By a theorem of Russo \cite{Russo} (see also \cite{PS}) applied to the triangular
lattice, this implies that the mean cluster size of the $+1$ and $-1$ clusters in,
say, ${\cal B}$ at odd $n \geq 1$ diverges.
It follows that, for any odd $n \geq 1$ and $x \in {\cal B}$
\begin{equation} \label{infty_B}
E(|C_x (n)|) \geq E(|C_x^{\cal B} (n)|) = \infty.
\end{equation}
Now, if $|C_x(n)| \neq 1$, then $C_x(n) = C_y(n)$ for some $y \in {\cal N}(x)$
and so
\begin{equation}
|C_x(n)| \leq \sum_{y \in {\cal N}(x)} |C_y(n)|.
\end{equation}
Taking expectations and noting that each neighbor $y$ of $x \in {\cal B}$ is in
${\cal A}$, we see that for odd $n \geq 1$, (\ref{infty_B}) for $x \in {\cal B}$ implies (\ref{inf_mean}) for all $x \in {\Bbb H}$.
The proof of (\ref{inf_mean}) for even $n \geq 2$ is similar with ${\cal A}$
and ${\cal B}$ interchanging their roles. \fbox{} \\

\bigskip
\bigskip

\section{Disordered ferromagnet}

In this section we study a disordered ferromagnet with i.i.d. couplings $J_{x,y}$
uniformly distributed between $0$ and $1$.
We will also have an external magnetic field $h>0$, so that the (formal) 
Hamiltonian is
\begin{equation}
{\cal H} = - \sum_{\langle x,y \rangle} J_{x,y} \sigma_x \sigma_y 
- h \sum_x \sigma_x.
\end{equation}
The initial density of $+1$ sites is $\lambda \in [0,1]$.
Let $P = P_{\sigma^0} \times P_{\omega} \times P_{\cal J}$ denote the joint
distribution on the space of $\sigma^0$'s, $\omega$'s and $ {\cal J}$'s.

\begin{remark}
The theorems and proofs of this section are valid for more general choices of the common
distribution of the $J_{x,y}$'s than uniform on $[0,1]$.
For example, Proposition \ref{h_geq_2} is valid for any distribution on $[0,1)$ and
Theorem \ref{dis_exp} is valid for any continuous distribution on $(0,1)$ that is symmetric about $1/2$.
\end{remark}

\begin{prop} \label{h_geq_2}
In the context of the disordered ferromagnet just described, if $h \geq 2$ and $\lambda \in [0,1]$,
almost surely, $\lim_{t \to \infty} \sigma^t_x = +1$ for all x.
Moreover, fixation happens exponentially fast in the following sense:
\begin{equation} \label{exp_fast}
\tilde{p}(t) \leq e^{- c t}
\end{equation}
for some $c \in (0, \infty)$.
\end{prop}

{\bf Proof.} Let $\hat{p}(t)$ denote the $P_{\sigma^0} \times P_{\omega}
\times P_{\cal J}$
probability that the origin is $-1$ at any time after $t$.
We claim that to prove both conclusions of the
proposition, it suffices
to show that $\hat{p}(t) \leq e^{-c t}$.
To see this, first note that  $\hat{p}(t) \to 0$ 
implies that the origin is eventually $+1$ almost surely, and
so by translation invariance, the same is true for any site; 
then note that
$\tilde{p}(t)
\leq \hat{p}(t)$ (since if a site flips then it must be minus
either just before or just after the flip).

Now notice that there is a positive density $\rho$ of sites $x$
such that $\sum_{y \in {\cal N} (x)} J_{x,y}$ $ < h$.
If the origin is one of those sites, then either it is $+1$ and will remain such,
or it is $-1$ and will flip to $+1$ the first time its clock rings, which happens
exponentially fast.
If the origin is not one of those sites, there is a closest site $X$ such that
$\sum_{y \in {\cal N} (X)} J_{X,y} < h$ and, since $\rho>0$, 
\begin{equation} \label{distance}
P_{\cal J} (||X|| > n) \leq e^{-c' n}
\end{equation}
for some $c' \in (0, \infty)$, where $|| X ||$
denotes the number of steps between the origin and site $X$ along a shortest path.
Next, notice that if a site $z$ has at least one neighbor
$y$ with $\sigma_y=+1$ when the clock at $z$ rings, $\sigma_z$ will flip to $+1$
if it is $-1$ and will remain $+1$ otherwise.
Then we have, for any $\alpha > 0$,
\begin{eqnarray} 
\hat{p}(t) & \leq & P_{\cal J} (||X|| > [\alpha t]) + 
\max_{0 \leq n \leq [\alpha t]} {\rm Prob} (T_0 + T_1 + \ldots + T_n > t) \nonumber \\
 & \leq & P_{\cal J} (||X|| > [\alpha t]) + 
{\rm Prob} (T_0 + T_1 + \ldots + T_{[\alpha t]} > t)
\end{eqnarray}
where $T_0, T_1, \ldots$ are independent exponential (mean one) random variables.
%
%
By choosing $\alpha$ small, both terms above can be seen to be
exponentially small in $t$, the first one using (\ref{distance}), the second one by
standard large deviation arguments. \fbox{} \\


When $0 \leq h<2$, it is easy to see that $\sigma^{\infty} = \lim_{t \to \infty} \sigma^t$
exists almost surely, but is not constant, that is, 
for any $\lambda \in (0,1)$ there is a positive
density in $\sigma^{\infty}$ of both $+1$ and $-1$ sites.
In this case, we can prove (in the next theorem) exponential fixation (in the above sense) only for
$h>1.5$, although we conjecture that exponential fixation happens for
all values of $h$ down to $h=0$.

\begin{thm} \label{dis_exp}
In the disordered ferromagnet with $\lambda \in [0,1]$, if $h > 1.5$, 
\begin{equation}
\tilde{p}(t) \leq e^{- c t}
\end{equation}
for some $c \in (0, \infty)$.
\end{thm}

To prove Theorem \ref{dis_exp}, we need the following definitions and lemmas:

\begin{df} \label{rings}
For each site $x$, let $\tau(x) = \{ \tau_n(x) : n \in {\Bbb N} \}$
denote the arrival times (i.e., times of clock rings) of the Poisson clock associated
with site $x$, arranged so that $\tau_k(x) < \tau_{k+1}(x)$ for each $k$.
Let $\tau^+ (x,t) = \inf \{ \tau_k(x) : \tau_k(x) > t \}$ be the next clock ring
at $x$ after time $t$.
\end{df}

\begin{df} \label{cascade} \emph{\cite{Howard}}
We call a sequence of sites $(x_1, \dots, x_n)$ with  $x_j \; $ $\in $ $\; {\cal N}(x_{j-1})$ for
$j = 2, \dots, n$ a \emph{plus-cascade} (of length $n$)
if, for a sequence of times $t_1, \dots, t_n$ defined by $t_1 = \tau_1(x_1)$
and, for $k>1$, $t_k = \tau^+ (x_k,t_{k-1})$, we have: $\sigma^t_{x_k}$ flips from $-1$
to $+1$ at time $t=t_k$ and, for $k<n$, does not flip again until possibly after time
$t_{k+1}$, and for $k>1$, $t_k>\tau_1(x_k)$.
\end{df}

Similarly, we define minus-cascades.
A sequence of neighboring sites $(x_1, \dots,$ $ x_n)$ is a cascade if it is either
a plus- or a minus-cascade.
If site $x$ does not flip at $\tau_1(x)$ ($x$'s first clock ring), then no
cascade begins at $x$ and we say that $x$ has an \emph{empty cascade}.
Notice that in the disordered models, cascades can split and merge.
As an immediate consequence of the definition, a
plus-cascade $(x_1, \dots, x_n)$ has the two following properties, which are
useful to keep in mind:
\begin{itemize}
\item for $1<k<n$, $x_k$ and $x_{k+1}$ are nearest neighbors,
\item for $1<k<n$, at the time $t = t_k$ when $\sigma_{x_k}$ flips to $+1$, 
$\sigma_{x_{k-1}}^t = +1$ and $\sigma_{x_{k+1}}^t = -1$.
\end{itemize}
Similar properties are valid for a minus-cascade.

\begin{lemma} \emph{\cite{Howard}}
Every flip belongs to some cascade. 
\end{lemma}

{\bf Proof.}
Consider a spin flip from (say) $-1$ to $+1$ at site $x$ at time $t$.
If $t = \tau_1(x)$, then the flip we are considering belongs to a cascade
starting at $x$.
If $t > \tau_1(x)$, then $t = \tau_k(x)$ for some $k \geq 2$.
In that case $\sigma_x^{t+} = - \sigma_x^{\tau_{(k-1)}(x)+}$, which implies that at least
one site $y_1 \in {\cal N}(x)$ experienced a spin flip from $-1$ to $+1$ at some time
$t_1 \in (\tau_{(k-1)}(x),t)$ and then did not flip again before time $t$.
We repeat the same procedure for $y_1$ and construct inductively a sequence of sites
$(y_0=x, y_1, \dots, y_n)$ and a sequence of times $t_0=t > t_1 > \dots > t_n$.
The procedure stops when $t_n = \tau_1(y_n)$.
It should be clear from the construction that $(y_n, \dots, y_1, y_0=x)$ is a cascade
for the sequence of times $t_n, \dots, t_1, t_0=t$. \fbox{} \\

\begin{lemma} \label{max_flips}
In the setting of Theorem \ref{dis_exp}, each site can flip at most $7$ times.
\end{lemma}

{\bf Proof.} By the translation invariance of the model, it is enough to prove
the claim for the origin $0$.
In order to do that, we will show that a minus-cascade can have at most length $2$.
Consider a minus-cascade whose first two sites are respectively $x_0$ and $y_0$.
Call $x_1$ and $x_2$ the other two neighboring sites of $x_0$, and $y_1$ and $y_2$
those of $y_0$.
When $x_0$ flips for the first time, from $+1$ to $-1$, by the definition of a
minus-cascade, $\sigma_{y_0}^{\tau_1(x_0)} = +1$, because $y_0$ is the
second site of the cascade.
Then, because of the strength of the field ($h>1.5$), for $\sigma_{x_0}$ to 
flip from $+1$ to $-1$ at time $t_1 = \tau_1(x_0)$, the two following conditions are
necessary:
$\sigma_{x_1}^{\tau_1(x_0)} = \sigma_{x_2}^{\tau_1(x_0)} = -1$,
and $J_{x_0,y_0} < 0.5$.
Now, for $\sigma_{y_0}$ to flip from $+1$ to $-1$ at time $t_2$, since $J_{x_0,y_0} < 0.5$,
it must be the case that $\sigma_{y_1} = \sigma_{y_2} = -1$.
Therefore, the cascade cannot proceed further, since when $\sigma_{y_0}$ flips
at time $t_2$, all its neighbors are already $-1$. 

Consider now the origin $0$.
It can flip from $+1$ to $-1$ at most three times, due to a cascade starting at
the origin itself plus at most two cascades starting at two of its neighbors or else
due to cascades starting at its three neighboring sites
(if there is a minus-cascade starting at the origin and then $\sigma_0$ flips back
to $+1$, that flip belongs to a plus-cascade involving some $y \in {\cal N}(0)$
and $y$ cannot then be the start of a minus-cascade reaching $0$).
In fact, a minus-cascade starting further away would not reach the origin.
Therefore, the origin can flip at most $7$ times. \fbox{} \\

{\bf Proof of Theorem \ref{dis_exp}}. Once again, we partition
the hexagonal lattice into the two triangular sublattices ${\cal A}$
and ${\cal B}$ as in the proof of Theorem \ref{exp_conv1}. 
In addition to the Poisson clock ${\bf C}_x$ of rate $1$ at site $x \in {\cal B}$, we also
assign three more clocks ${\bf C}_{x,y}$ with rate $1/3$ each and associated respectively
with the neighboring sites $y \in {\cal N} (x)$, all of them being in ${\cal A}$.
(We remark that these clocks ${\bf C}_{x,y}$ are different from the ones used in the proof
of Theorem \ref{exp_conv1}.)
All the clocks ${\bf C}_x, {\bf C}_{x,y}$ for $x \in {\cal B}$ and $y \in {\cal N} (x)$ are
independent.
We now define a new dynamics according to the following rules:
\begin{itemize}
\item A site $x \in {\cal B}$ attempts a spin flip when ${\bf C}_x$ rings.
\item A site $y \in {\cal A}$ attempts a spin flip whenever any one of
the clocks ${\bf C}_{x,y}$  assigned to one of its neighbors $x \in {\cal N} (y)$
(in sublattice ${\cal B}$) rings.
\item The rules for accepting a spin flip are the same as before.
\end{itemize}

Since the clocks ${\bf C}_{x,y}$ used for sites $y$ in ${\cal A}$ have each rate $1/3$,
the resulting dynamics is the same as before.
We fix a time $T>0$ and say that $x \in {\cal B}$ is a \emph{good} site if the two following
conditions hold:
\begin{enumerate}
\item $\sum_{y \in {\cal N} (x)} J_{x,y} < h$, and
\item during the time interval $[0, T]$, ${\bf C}_x$ rings and following that
all three clocks ${\bf C}_{x,y}$ for $y \in {\cal N} (x)$ also ring during $[0, T]$.
\end{enumerate}
A site that is not good will be called \emph{bad}.
Clearly, because of condition 1 and the external magnetic field, 
good sites that are ever $+1$ are then stable for the dynamics.
Condition 2 implies that a good site $x \in {\cal B}$ that is $-1$ at time zero will have
flipped to become a stable $+1$ by time $T$, and moreover it implies that a site 
$y \in {\cal A}$ with two
good sites as ${\Bbb H}$-neighbors will be a stable $+1$ by time $T$, regardless of its
value at time zero.

The event corresponding to condition 1 has probability strictly larger than $1/2$ 
because of the distribution of the couplings and our assumption that
$h>1.5$.
The event corresponding to condition 2 has probability approaching one as $T \to \infty$.
So, choosing $T$ large enough, we will have that the probability that a site $x$
is good is strictly larger than $1/2$.
Thus good sites percolate in the triangular lattice ${\cal B}$ and therefore, for $x \in {\cal B}$
(say, the origin), if $B_x^{\cal B}$ is a ${\cal B}$-cluster of bad sites,
\begin{equation}
(P_{\omega} \times P_{\cal J}) (|B_x^{\cal B}| > n) \leq e^{- c'' n}
\end{equation}
for some $c'' \in (0, \infty)$.
Moreover, the boundary $\partial B_x^{\cal B}$ of  the cluster $B_x^{\cal B}$ is the union
of (one or more) ${\cal B}$-connected loops of good sites.
Let $\partial_x^{\cal A} = \{ y \in {\cal A} : \exists \, \zeta_1, \zeta_2 \in 
\partial B_x^{\cal B} \, \text{s.t.} \, ||\zeta_1 - y|| = ||\zeta_2 - y|| = 1 \}$
and define $\partial_x = \partial B_x^{\cal B} \cup \partial_x^{\cal A}$.
Then $\partial_x$ is  the union of (one or more) ${\Bbb H}$-connected sets and it
completely surrounds the ${\Bbb H}$-connected set $B_x = B_x^{\cal B} \cup B_x^{\cal A}$, 
where $B_x^{\cal A} = \{ y \in {\cal A} : \exists \, \zeta_1, \zeta_2 \in 
B_x^{\cal B} \, \text{s.t.} \, ||\zeta_1 - y|| = ||\zeta_2 - y|| = 1 \}$.
Clearly $|B_x| \leq 3 \, |B_x^{\cal B}|$ and therefore
\begin{equation}
(P_{\omega} \times P_{\cal J}) (|B_x| > n) \leq e^{- c'' n/3} = e^{- c' n}
\end{equation}
for some $c' \in (0, \infty)$.

Because of condition 2 above, by time $T$ all the sites in $\partial_x$ have fixated:
the sites in $\partial B_x^{\cal B}$ because they are good sites, those in $\partial_x^{\cal A}$
because they have two good sites as ${\Bbb H}$-neighbors.
To conclude the proof, observe that if the origin has not fixated by time $T$, it must
belong to a connected set $B_0$ surrounded by fixated sites and with 
\begin{equation} \label{dis_expon}
(P_{\omega} \times P_{\cal J}) (|B_0| > n) \leq e^{- c' n}.
\end{equation}
%
%
By Lemma \ref{max_flips}, the maximum number of spin flips allowed inside $B_0$ is
$7 \, |B_0|$.
Moreover, as long as a single unstable site is present, it will flip at its next
clock ring.
It follows that an upper bound for the time for the sites in $B_0$ to fixate is
given by a sum of $7 \, |B_0|$ exponential (mean one) random variables $T_1, T_2, \dots$ .
More precisely, for any $\alpha >0$ and $t>T$,
\begin{eqnarray}
\tilde{p}(t) & \leq & (P_{\omega} \times P_{\cal J}) (|B_0| > [\alpha \, (t-T)]/7)\nonumber \\
&\;  & \qquad + 
\max_{1 \leq n \leq [\alpha \, (t-T)]} {\rm Prob}(T_1 + \ldots + T_n > (t-T)) \nonumber \\
& = & (P_{\omega} \times P_{\cal J}) (|B_0| > [\alpha \, (t-T)]/7) \nonumber \\
&\;  & \qquad+ {\rm Prob}(T_1 + \ldots +
T_{[\alpha \, (t-T)]} > (t-T)). \label{dis_p_t} 
\end{eqnarray}
By choosing $0 < \alpha < 1$, both terms in (\ref{dis_p_t}) can be seen to be
exponentially small in $t$, the first one using (\ref{dis_expon}), the second one by
standard large deviation arguments. \fbox{} \\

\bigskip
\noindent  {\bf Acknowledgments.} Research partially supported by the U.S. NSF under grants
DMS-98-02310 and DMS-01-02587 (F. Camia),  DMS-98-03267  and DMS-01-04278 (C.M. Newman),  Faperj grant E-26/151.905/2000,
Pronex and CNPq (V. Sidoravicius). A portion of this research was done while some of us were visitors at ETHZ, the Courant
Institute and IMPA; we thank Alain Sznitman and these institutions for their hospitality.

%
%

\end{document}